\documentclass[11pt,amsfonts]{article}
\usepackage{graphicx,amssymb,eucal,mathrsfs,url}
\usepackage{latexsym,amsmath,amscd,epsfig,xcolor}
\usepackage[all]{xy}
\usepackage{ntheorem,tikz}
\usepackage{epsf,graphicx,pgf,etex,amscd,pstricks}

\usepackage[hypertexnames=false,backref=page,pdftex,
 	pdfpagemode=UseNone,
 	breaklinks=true,
 	extension=pdf,
 	colorlinks=true,
 	linkcolor=blue,
 	citecolor=red,
 	urlcolor=blue,
 ]{hyperref}

\newcommand{\Z}{\mathbb{Z}}
\newcommand{\Q}{\mathbb{Q}}
\newcommand{\R}{\mathbb{R}}
\newcommand{\C}{\mathbb{C}}

\newtheorem{main}{Theorem}

\newtheorem{question}[main]{Question}
\newtheorem{prop}[main]{Proposition}
\newtheorem{lemma}[main]{Lemma}
\newtheorem{cor}[main]{Corollary}

\theorembodyfont{\upshape}

\newtheorem{defi}[main]{Definition}
\newtheorem{rem}[main]{Remark}

\renewcommand{\thefootnote}{\alph{footnote}}

\title{Subgroups of hyperbolic groups, finiteness properties and complex hyperbolic lattices}
\date{April 2022}
\author{Claudio Llosa Isenrich$^{{\rm a}}$ and Pierre Py}

\begin{document}

\maketitle

\begin{abstract} We prove that in a cocompact complex hyperbolic arithmetic lattice $\Gamma < {\rm PU}(m,1)$ of the simplest type, deep enough finite index subgroups admit plenty of homomorphisms to $\Z$ with kernel of type $\mathscr{F}_{m-1}$ but not of type $\mathscr{F}_{m}$. This provides many finitely presented non-hyperbolic subgroups of hyperbolic groups and answers an old question of Brady. Our method also yields a proof of a special case of Singer's conjecture for aspherical K\"ahler manifolds. 
\end{abstract}

\footnotetext[1]{The author gratefully acknowledges funding by the DFG 281869850 (RTG 2229).}

\renewcommand{\thefootnote}{\arabic{footnote}}


\section{Introduction}\label{sec:Intro}

A classifying space for a group $G$, or $K(G,1)$, is an aspherical ${\rm CW}$-complex with fundamental group $G$. Following Wall~\cite{wall}, we say that $G$ is of type $\mathscr{F}_{n}$ if it has a $K(G,1)$ with finite $n$-skeleton. One usually refers to property $\mathscr{F}_{n}$ as a {\it finiteness property} for the group $G$. Property $\mathscr{F}_{1}$ (resp. $\mathscr{F}_{2}$) is equivalent to being finitely generated (resp. finitely presented). We say that $G$ is of type $\mathscr{F}_{\infty}$ if it is $\mathscr{F}_n$ for all $n$ and that $G$ is of type $\mathscr{F}$ if it admits a finite $K(G,1)$. For each integer $n$, there are groups of type $\mathscr{F}_{n}$ but not of type $\mathscr{F}_{n+1}$ \cite{bieri-76,stallings}. There are also other families of finiteness properties for groups~\cite{bb,brown}; the only one which we shall refer to below is property ${\rm FP}_{n}(\Q)$, see~\cite{bb} for its definition. If a group is of type $\mathscr{F}_{n}$, then it is of type ${\rm FP}_{n}(\Q)$, but the converse implication does not hold. 

Using methods from {\it complex geometry}, we prove for each integer $n\geq 1$ the existence of Gromov hyperbolic groups containing coabelian subgroups of type $\mathscr{F}_{n}$ but not of type $\mathscr{F}_{n+1}$. We now state our main results more precisely, before explaining the historical context and motivation behind them. 

The {\it homotopical BNSR invariants} of a finitely generated group $G$ form a sequence $(\Sigma^{j}(G))_{j\ge 1}$ of open subsets of the {\it character sphere} $$S(G):=H^{1}(G,\R)-\{0\}/\R_{+}^{\ast}.$$ They were introduced by Bieri, Neumann and Strebel in~\cite{bns} (for $j=1$) and by Renz in~\cite{Renz-thesis} (for $j\ge 2$). Their definition is recalled in Section~\ref{sec:bnsr}. We simply mention here that they encode in some sense the finiteness properties of kernels of homomorphisms from $G$ to $\R$. In what follows we write $[\xi]$ for the image in the sphere $S(G)$ of a nonzero class $\xi \in H^{1}(G,\R)$. 

\begin{main}\label{bigsigmahd} Let $m\ge 2$ and let $\Gamma < {\rm PU}(m,1)$ be a torsion-free cocompact arithmetic lattice of the simplest type. Then $\Gamma$ has a finite index subgroup $\Gamma_0$ with the following property. For every finite index subgroup $\Gamma_1 < \Gamma_0$ the $(m-1)$-th BNSR invariant $\Sigma^{m-1}(\Gamma_{1})$ of $\Gamma_{1}$ is dense in the character sphere $S(\Gamma_{1})$. In particular, every rational class $\xi \in H^{1}(\Gamma_{1},\mathbb{Q})$ such that $[\xi ]$ is contained in the dense open set $\Sigma^{m-1}(\Gamma_{1})\cap - \Sigma^{m-1}(\Gamma_{1})$ satisfies that ${\rm ker}(\xi) <\Gamma_{1}$ is a group of type $\mathscr{F}_{m-1}$ but not of type ${\rm FP}_{m}(\Q)$.   
\end{main}

{\it Arithmetic lattices of the simplest type} are the lattices associated to Hermitian forms with coefficients in a purely imaginary quadratic extension of a totally real number field. See Section~\ref{sec:defsimplesttype} for their definition. Before going further we also introduce the following classical:

\begin{defi}\label{def:cor} Two groups are said to be commensurable if they have isomorphic finite index subgroups.
\end{defi}

Since arithmetic lattices of the simplest type in ${\rm PU}(m,1)$ ($m\ge 2$) form an infinite family of commensurability classes (see Section~\ref{sec:defsimplesttype}), Theorem~\ref{bigsigmahd} has the following consequence. 

\begin{cor}\label{corollairerationalclass} Let $n\ge 2$ be an integer. There exist infinitely many hyperbolic groups $(G_{n,j})_{j\ge 0}$ and homomorphisms $\phi_{n,j} : G_{n,j}\to \mathbb{Z}$ such that the $G_{n,j}$'s are pairwise noncommensurable and such that the kernel of $\phi_{n,j}$ is of type $\mathscr{F}_{n}$ but not of type ${\rm FP}_{n+1}(\Q)$. 
\end{cor}

We now move on to a more detailed introduction. A finitely generated group is called {\emph{hyperbolic}} if its Cayley graph with respect to some finite generating set is $\delta$-hyperbolic~\cite{cdp,ghysharpe}. Introduced in Gromov's seminal essay~\cite{gromov}, the class of hyperbolic groups has attracted much attention. Indeed these groups satisfy many nice properties. For instance, they have solvable word and conjugacy problem, they do not contain any $\Z^2$-subgroups, they satisfy the Tits alternative and every hyperbolic group is of type $\mathscr{F}_{\infty}$ (and of type $\mathscr{F}$ in the torsion-free case). It is natural to look for properties of hyperbolic groups which are satisfied by all of their subgroups. While some of the aforementioned properties obviously pass to all subgroups, answering this question for others turns out to be surprisingly difficult. It is in this general context that Brady asked in 1999 the following: 
\begin{question}[{Brady \cite{brady}}]\label{qu:brady}
 For a given integer $n\geq 3$, does there exist a hyperbolic group $G$ which has a subgroup $H<G$ of type $\mathscr{F}_n$, but not of type $\mathscr{F}_{n+1}$?
\end{question}

Preceding Brady's question, in 1982 Rips proved the existence of {\it non-coherent} hyperbolic groups \cite{rips}, i.e. the existence of subgroups of hyperbolic groups of type $\mathscr{F}_1$ and not $\mathscr{F}_2$; many more examples have been constructed since by different authors. Gersten proved that every finitely presented subgroup of a hyperbolic group of cohomological dimension two is hyperbolic \cite{gersten}. However, this result is not true for subgroups of hyperbolic groups of higher cohomological dimension. The existence of finitely presented non-hyperbolic subgroups of hyperbolic groups was first suggested by Gromov in his 1987 essay \cite{gromov}. In~\cite{brady}, Brady constructed the first example of such a subgroup. It is of type $\mathscr{F}_2$ but not of type $\mathscr{F}_3$ (thus answering the question above for $n=2$) and arises as the fundamental group of a generic fibre of a map from a ramified covering of a direct product of three graphs onto the circle. More examples of subgroups of hyperbolic groups of type $\mathscr{F}_2$ and not $\mathscr{F}_3$ have been constructed since by Lodha \cite{lodha} and Kropholler \cite{krophollergardam}. 

Very recently the authors of the present article, together with Martelli~\cite{limp} showed the existence of subgroups of hyperbolic groups of type $\mathscr{F}_3$ and not $\mathscr{F}_4$. These examples are obtained by starting from a cusped real hyperbolic 8-manifold $M^8$ and a map $f:M^8\to S^1$ with $\ker(f_{\ast}:\pi_1(M^8)\to \Z)$ of type $\mathscr{F}_3$ and not $\mathscr{F}_4$ and then considering certain Dehn fillings of $M^8$. The map $f$ was constructed by Italiano, Martelli and Migliorini \cite{itmami}. Before~\cite{limp}, Italiano, Martelli and Migliorini had already built a fibration from a cusped hyperbolic 5-manifold to the circle and from it had produced a non-hyperbolic finitely presented subgroup of type $\mathscr{F}$ of a hyperbolic group \cite{itmami2}. This constituted fundamental progress in the area.

All of the aforementioned constructions rely on Bestvina--Brady Morse theory~\cite{bb}, which requires checking some combinatorial conditions. While one may hope that similar methods will allow one to answer Brady's question for $n\geq 4$ using e.g. cubulated lattices in ${\rm PO}(2n,1)$ or right-angled Coxeter groups, we follow a different path here. We use complex hyperbolic lattices instead of real hyperbolic lattices and apply complex Morse theory (also known as {\it Lefschetz theory}) instead of Bestvina--Brady Morse theory. This allows us to answer elegantly Brady's question for all $n$. In this context, the high connectivity of the kernels that we study appears as a natural consequence of Lefschetz theory. For earlier uses of this theory to study finiteness properties of groups, see~\cite{dps,kapovich,Llo-17,nicolaspy}.   

One may raise analogues of Question~\ref{qu:brady} for other finiteness properties of groups. Motivated by this, we mention a relation between our work and the work of Fisher and Kielak~\cite{fisher,kielak,kielakoberwolfach}. In what follows we will always denote by $b_{i}^{(2)}(G)$ the $i$-th $\ell^{2}$-Betti number of a group $G$; see~\cite{chegro1986,luck} for the definition. A classical theorem of L\"uck~\cite[Th. 3.3 (4)]{luck1998-II} implies that if $G$ is any group, if $\phi :  G \to \Z^{k\ge 1}$ is a surjective homomorphism and if $b_{i}^{(2)}(G)\neq 0$ then the kernel of $\phi$ cannot be of type ${\rm FP}_{i}(\Q)$. See e.g.~\cite{fisher} or~\cite[Prop. 14]{limp} for a proof of this fact. Hence the presence of nonvanishing $\ell^{2}$-Betti numbers provides an ``upper bound"  for the finiteness properties of coabelian normal subgroups. For RFRS groups this bound is sharp in the following sense:

\begin{main}\label{th:kf} (Fisher, Kielak) Let $G$ be a group which is virtually RFRS and of type ${\rm FP}_{n}(\Q)$. Then the following conditions are equivalent:
\begin{enumerate} 
\item there exists a finite index subgroup $G_1 < G$ and a homomorphism $\phi : G_{1} \to \Z$ whose kernel is of type ${\rm FP}_{n}(\Q)$;
\item the $\ell^{2}$-Betti number $b_{i}^{(2)}(G)$ vanishes for $0\le i \le n$.  
\end{enumerate} 
\end{main}

We refer to~\cite{agol} for the definition of RFRS groups. Kielak~\cite{kielak} proved the case $i=1$ of the theorem above and conjectured the higher degree version~\cite{kielakoberwolfach}, which was subsequently proved by Fisher~\cite{fisher}. Theorem~\ref{bigsigmahd} is philosophically similar to Theorem~\ref{th:kf}. Although we deal with a specific class of groups (cocompact complex hyperbolic arithmetic lattices of the simplest type), which are not known to be RFRS, the conclusion is close in spirit. We prove that sufficiently deep finite index subgroups in such lattices admit homomorphisms to $\Z$ (in fact plenty of them) whose finiteness properties are as ``good" as allowed by the $\ell^{2}$-Betti numbers. Indeed, it is known that the $\ell^{2}$-Betti numbers $b_{i}^{(2)}(\Gamma)$ of a cocompact lattice $\Gamma < {\rm PU}(m,1)$ are all zero except for the $m$-th one which is nonzero~\cite{borel}. The difference between the two theorems is that we use the stronger property $\mathscr{F}_{n}$ instead of property ${\rm FP}_{n}(\Q)$. Note also that in Theorem~\ref{bigsigmahd}, the fact that the kernels of the rational characters lying in the set $\Sigma^{m-1}(\Gamma_{1})\cap -\Sigma^{m-1}(\Gamma_{1})$ are not of type ${\rm FP}_{m}(\Q)$ can alternatively be shown by a topological argument not relying on $\ell^{2}$-Betti numbers, see~\cite[Prop. 21]{limp}. Beyond the philosophical similarity between Theorems~\ref{bigsigmahd} and ~\ref{th:kf} we would also like to point out that one of the consequences of Theorem ~\ref{th:kf} is the existence of subgroups of hyperbolic groups of type ${\rm FP}_n(\Q)$, but not ${\rm FP}_{n+1}(\Q)$ for all integers $n\geq 1$ \cite[Prop. 19]{limp}.

We now move to the complex geometric setting. We shall deduce Theorem~\ref{bigsigmahd} from the following result, dealing with arbitrary closed aspherical K\"ahler manifolds.

\begin{main}\label{theo:varcomplexabs} Let $X$ be a closed aspherical K\"ahler manifold of complex dimension $m\ge 2$. 
\begin{enumerate}
\item Let $\beta$ be a holomorphic $1$-form on $X$ with finitely many zeros. Then the cohomology class $b=[Re(\beta)] \in H^{1}(X,\R)\simeq H^{1}(\pi_{1}(X),\R)$ lies in $\Sigma^{m-1}(\pi_{1}(X))\cap -\Sigma^{m-1}(\pi_{1}(X))$. If $b$ is rational then its kernel is of type $\mathscr{F}_{m-1}$; if furthermore the Euler characteristic of $X$ is nonzero, the kernel of $b$ is not of type ${\rm FP}_{m}(\Q)$. 
\item Let $\psi : X \to A$ be a holomorphic map to a complex torus. Assume that $\psi$ is a finite map. Let $\alpha$ be a holomorphic $1$-form on $A$ which does not vanish on any nontrivial subtorus of $A$. Then the form $\psi^{\ast}\alpha$ has finitely many zeros. Consequently, the class $[Re (\psi^{\ast}\alpha)]$ lies in $\Sigma^{m-1}(\pi_{1}(X))\cap -\Sigma^{m-1}(\pi_{1}(X))$. 
\end{enumerate} 
\end{main}

We recall that a map $f : X \to Y$ between two topological spaces is said to be finite if each of its fibers $f^{-1}(y)$ ($y\in Y$) is a finite set.

The proof of Theorem~\ref{theo:varcomplexabs} relies on the work of Simpson~\cite{simpson}. Note that Delzant~\cite{delzant} used the same work of Simpson to give a complete description of the BNS invariant $\Sigma^{1}$ for K\"ahler groups. Our results can be seen as a higher degree generalization of Delzant's work. Of course we have to make additional hypotheses (we consider an aspherical K\"ahler manifold and finite maps to complex tori) and we do not describe completely the invariant $\Sigma^{m-1}$ but just exhibit a large set contained in it. 

Throughout this text, we will identify the group ${\rm PU}(m,1)$ with the group of holomorphic automorphisms of the unit ball $B$ of $\C^{m}$. If $\Gamma < {\rm PU}(m,1)$ is a lattice, we will denote by $X_{\Gamma}$ the corresponding quotient:
\begin{equation}
X_{\Gamma}:=B/\Gamma .    
\end{equation}
The link between Theorem~\ref{theo:varcomplexabs} and Theorem~\ref{bigsigmahd} is provided by the following two facts. Firstly, for a lattice $\Gamma < {\rm PU}(m,1)$ satisfying the hypotheses of Theorem~\ref{bigsigmahd}, and for a deep enough finite index subgroup $\Gamma_0 < \Gamma$, the first Betti number of $\Gamma_0$ is positive~\cite{kazhdan} and the Albanese map of $X_{\Gamma_{0}}$ is an immersion, thus a finite map. The latter fact is due to Eyssidieux~\cite{eyssidieux} and we give a short account of it in Section~\ref{sec:chlattices}. Secondly, keeping the notations from Theorem~\ref{theo:varcomplexabs} and assuming that $A$ is the Albanese torus of $X$, the condition ``$\alpha$ does not vanish on any nontrivial subtorus" is a condition satisfied by a dense set of classes $$a= [Re (\psi^{\ast}\alpha)]\in H^{1}(X,\R),$$ see Proposition~\ref{prop:lacestod}. For the reader more familiar with hyperbolic groups than with K\"ahler manifolds, we mention that the Albanese torus ${\rm Alb}(X)$ of a closed K\"ahler manifold $X$ is a complex torus of complex dimension $\frac{1}{2}b_{1}(X)$ which comes with a natural holomorphic map ${\rm alb}_{X} : X \to {\rm Alb}(X)$ inducing an isomorphism between the first real (co)homology groups of $X$ and ${\rm Alb}(X)$. We will recall its definition in Section~\ref{sec:chlattices}.

Our results raise the question of whether one can give a complete description of the invariants $\Sigma^{j}$ ($j\ge 2$) for fundamental groups of closed aspherical K\"ahler manifolds, or at least in some specific situations. This problem is related to a question of Kotschick~\cite[Question 15]{kotschick}. Note also that Friedl and Vidussi conjecture that $\Sigma^{2}(\pi_{1}(X))$ is always empty when $X$ is a closed aspherical K\"ahler surface of nonzero Euler characteristic~\cite[p. 53]{friedlvidussi}. So far we have seen that a result of L\"uck on $\ell^{2}$-Betti numbers gives restrictions on finiteness properties of kernels of homomorphisms to Abelian groups. This result implies that if $b_{i}^{(2)}(G)$ is nonzero for a group $G$ which is of type $\mathscr{F}_{i}$, then $\Sigma^{i}(G)\cap -\Sigma^{i}(G)$ must be empty. In our setting, we ask:

\begin{question}
Let $X$ be a closed aspherical K\"ahler manifold of nonzero Euler characteristic. Let $m={\rm dim}_{\C}\, X$. Is it true that $\Sigma^{m}(\pi_{1}(X))$ is empty?
\end{question}
This question generalizes the conjecture by Friedl and Vidussi. In this direction we prove:

\begin{main}\label{th:fvsigmadeux} Let $X$ be a closed aspherical K\"ahler manifold of nonzero Euler characteristic. Let $m={\rm dim}_{\C} \, X$. Assume that the Albanese map of $X$ is finite. Then $\Sigma^{m}(\pi_{1}(X))$ is empty.
\end{main}

Combined with the properties of the Albanese map of arithmetic ball quotients alluded to above, Theorem~\ref{th:fvsigmadeux} has the following consequence. 

\begin{cor}\label{cor:sigmadeuxvide}
Let $\Gamma < {\rm PU}(m,1)$ be a torsion-free cocompact arithmetic lattice of the simplest type. Then $\Sigma^{m}(\Gamma)$ is empty. 
\end{cor}

Finally, we also observe that our results naturally imply a special case of the {\it Singer conjecture}~\cite{singer}, in the context of K\"ahler manifolds. This conjecture states that a closed aspherical $n$-dimensional manifold $M$ satisfies $b_{i}^{(2)}(M)=0$ if $2i$ is distinct from the dimension of $M$. See~\cite{luck} for a survey of known cases. In the context of K\"ahler manifolds, the most important result is Gromov's theorem~\cite{gromovkh} stating that a {\it K\"ahler hyperbolic} manifold satisfies the Singer conjecture (without the asphericity hypothesis). Consider now a closed K\"ahler manifold $X$ admitting a holomorphic $1$-form with finitely many zeros. This implies that the top Chern number of $T^{\ast}X$ is nonnegative, hence 
\begin{equation}\label{eq:hopfthurston}
(-1)^{m} \chi (X)\ge 0,
\end{equation}
where $m={\rm dim}_{\C} \, X$ and $\chi (X)$ is the Euler characteristic. In other words, $X$ satisfies the conclusion of the Hopf-Chern-Thurston conjecture (see~\cite[Ch. 16]{davisbook}). Recall that this conjecture states that Equation~\eqref{eq:hopfthurston} holds for all real aspherical $2m$-dimensional closed manifolds and is implied by Singer's conjecture. Assuming that $X$ is aspherical, our methods naturally yield the following stronger conclusion.   

\begin{main}\label{th:singerpartial} Let $X$ be a closed aspherical K\"ahler manifold. Assume that $X$ carries a holomorphic $1$-form $\alpha$ with finitely many zeros. Then $X$ satisfies the Singer conjecture. 
\end{main} 

The article is organized as follows. In Section~\ref{sec:bnsrds}, we recall the definition of the BNSR invariants of a group $G$ and then prove Theorems~\ref{theo:varcomplexabs}, \ref{th:fvsigmadeux} and \ref{th:singerpartial}. In Section~\ref{sec:chlattices}, we study the Albanese map of arithmetic quotients of complex hyperbolic space, then we recall the definition of arithmetic lattices of the simplest type, and finally we prove Theorem~\ref{bigsigmahd} and Corollaries~\ref{corollairerationalclass} and~\ref{cor:sigmadeuxvide}.


\section{Lefschetz theory and finiteness properties}\label{sec:bnsrds}

\subsection{The BNSR invariants}\label{sec:bnsr}

We recall here the definition of the BNSR invariants of a finitely generated group $G$~\cite{bns,Renz-thesis}. As in the introduction, we set
$$S(G):= H^{1}(G,\R)-\{0\} /\mathbb{R}_{+}^{\ast}$$
where $\R_{+}^{\ast}$ acts by scalar multiplication on $H^{1}(G,\R)$. The set $S(G)$ is called the {\it character sphere} of $G$. The BNSR invariants form a decreasing family of open subsets
$$\Sigma^{k}(G)\subset  \cdots \subset \Sigma^{2}(G)\subset \Sigma^{1}(G) \subset S(G)$$
of the character sphere of $G$. 
The invariant $\Sigma^{m}(G)$ is defined for groups of type $\mathscr{F}_{m}$ only. Here we will focus on groups which are fundamental groups of closed aspherical manifolds, hence this condition will be automatically satisfied. 

So let $M$ be a closed aspherical manifold with fundamental group $G$ and universal cover $\pi : \widehat{M} \to M$. Let $\chi : G \to \R$ be a nonzero character. We pick a closed $1$-form $u$ on $M$ representing $\chi$ and write $\pi^{\ast}u=df$ for some smooth function $f : \widehat{M} \to \R$. We let $\widehat{M}_{d}=f^{-1}([d,\infty))$. 

\begin{defi} We say that $\widehat{M}_{d}$ is essentially $m$-connected if there exists a real number $r\ge 0$ such that the inclusion map $\widehat{M}_{d}\to \widehat{M}_{d-r}$ induces the zero map $\pi_{i}(\widehat{M}_{d})\to \pi_{i}(\widehat{M}_{d-r})$ on homotopy groups for $i\le m$ (for $i=0$ this means by convention that the image of the map on $\pi_0$ is a singleton). 
\end{defi}

One can easily prove that $\widehat{M}_{d}$ is essentially $m$-connected for some real number $d$ if and only if it is essentially $m$-connected for all $d$. Moreover, the fact that $\widehat{M}_{d}$ is essentially $m$-connected only depends on the ray $[\chi]\in S(G)$, i.e. it does not depend on the choice of $u$ and $f$ and is unaffected if we multiply $\chi$ by a positive real number. We can thus introduce the following: 

\begin{defi}\label{def:bnsropla} The set $\Sigma^{m}(G)\subset S(G)$ consists of the rays $[\chi]$ such that $\widehat{M}_{d}$ is essentially $(m-1)$-connected for some $d\in \mathbb{R}$.
\end{defi}

When $m=1$, one can also define $\Sigma^{1}(G)$ by studying the connectivity of a certain subgraph of a Cayley graph of $G$, see e.g.~\cite[Def. 3.9]{kielak}. A proof of the fact that the definition given in~\cite{kielak} is equivalent to Definition~\ref{def:bnsropla} (for $m=1$) can be found in~\cite{bieristrebel}. In general, the definition of $\Sigma^{m}(G)$ is given using a ${\rm CW}$-complex $Q$ which is a $K(G,1)$ (instead of our manifold $M$), $Q$ being assumed to have finite $m$-skeleton, and working with a function $F$ from the universal cover of $Q$ to $\R$ which is $\chi$-equivariant (i.e. satisfies $F(g\cdot x)=F(x)+\chi (g)$). The definition is then formulated in terms of the homotopical properties of the sets $\{F \ge d\}$ as above, see e.g.~\cite[Remark 6.5]{BieriRenz} or~\cite{Renz-thesis}. Assuming that $G$ is the fundamental group of a closed aspherical manifold, one checks readily that our definition is equivalent to the original one \cite[Prop. B2.1]{bieristrebel}.   

We close this short presentation by mentioning two results about the invariants $\Sigma^{m}(G)$. First they are open subsets of the sphere $S(G)$, see~\cite{bieristrebel} and~\cite[\S IV.2]{Renz-thesis}. Second, we have the following~\cite[\S V.2]{Renz-thesis}:

\begin{main}\label{th:ratfn} Let $[\chi]\in S(G)$ be a rational point. Then the kernel of $\chi$ is of type $\mathscr{F}_{m}$ if and only if $[\chi]\in \Sigma^{m}(G)\cap -\Sigma^{m}(G)$. 
\end{main}

When dealing with certain rational cohomology classes below, we will propose two proofs of the fact that their kernels are of type $\mathscr{F}_{m}$. One is based on Theorem~\ref{th:ratfn}, the other is based on a direct complex Morse theory argument going back to~\cite{dps}.

\subsection{A higher dimensional version of Delzant's theorem}\label{hdvdt}

In this section we prove Theorem~\ref{theo:varcomplexabs}. For the moment we keep its notations and asumptions, except for one thing: we do not assume yet that $X$ is aspherical. We first deal with the second point. It reduces to the first, thanks to the following proposition, which already appears in~\cite{simpson}.

\begin{prop}\label{prop:alsago}
Under  the assumptions of Theorem~\ref{theo:varcomplexabs}, the form $\psi^{\ast}\alpha$ has only finitely many zeros on $X$.
\end{prop}

\noindent {\it Proof.} Let $Z$ be a connected component of the set of zeros of $\psi^{\ast}\alpha$. If $\psi (Z)$ is positive dimensional, it generates a nontrivial subtorus of $A$ on which $\alpha$ vanishes (see~\cite[VIII.1]{debarre}). This is a contradiction. Hence $\psi (Z)$ must be zero dimensional and thus a point. Since $\psi$ is finite, this implies that $Z$ is a point.\hfill $\Box$

We now deal with the first point of Theorem~\ref{theo:varcomplexabs}. Since $\beta$ and $-\beta$ both have finitely many zeros, it is enough to show that $b=[Re (\beta )]\in \Sigma^{m-1}(\pi_{1}(X))$. We then have the following result due to Simpson~\cite{simpson}. 

\begin{main}\label{th:simpsonpcisoles} Let $Y$ be a compact K\"ahler manifold of complex dimension $m\ge 2$ and $\beta$ be a holomorphic $1$-form on $Y$ with finitely many zeros. Let $\widehat{Y}$ be the universal cover of $Y$ and $f : \widehat{Y}\to \R$ be a primitive of the lift to $\widehat{Y}$ of the form $Re(\beta)$. Then for all real numbers $c, d$ such that $c\le d$ the inclusion 
$$f^{-1}([d,\infty))\subset f^{-1}([c,\infty))$$
induces an isomorphism on $\pi_{i}$ for $i\le m-2$ and a surjection on $\pi_{m-1}$. 
\end{main} 

This result is essentially contained in Theorem~17 and the subsequent remark from~\cite{simpson}. However, in~\cite{simpson} Simpson deals with a more general situation: he considers the homotopical properties of the level sets of the primitive of (the lift to $\widehat{X}$ of) either a harmonic form or a holomorphic form and he also allows the case where one works on a covering space of a quasi-projective variety (with extra assumptions). This makes the proof more involved. To make this text more self-contained, we shall give a quick proof of Theorem~\ref{th:simpsonpcisoles} in Section \ref{sec:simpson}, following Simpson's approach. 

We now explain how to conclude the proof of the first item of Theorem~\ref{theo:varcomplexabs}, using Theorem~\ref{th:simpsonpcisoles}. Assume that $X$ is aspherical. Let $\widehat{X}$ be the universal cover of $X$ and let $f : \widehat{X} \to \R$ be a primitive of the lift of $\beta$ to $\widehat{X}$. We set $\widehat{X}_{d}=f^{-1}([d,\infty))$, as in Section~\ref{sec:bnsr}. We shall prove that 
\begin{equation}
\pi_{i}(\widehat{X}_{d})=0
\end{equation}
for $i\le m-2$ and every real number $d$; this obviously implies that $\widehat{X}_{d}$ is essentially $(m-2)$-connected. Theorem~\ref{th:simpsonpcisoles} implies that $\widehat{X}_{d}$ is path-connected for all $d$. So let $1\le i\le m-2$ and $\xi : S^{i} \to \widehat{X}_{d}$ be a continuous map representing a class in $\pi_{i}(\widehat{X}_{d})$. Since $\widehat{X}$ is contractible, $\xi$ extends to a continuous map $\overline{\xi} : B^{i+1} \to \widehat{X}$. If $c:={\rm min}\{ d, {\rm inf}_{B^{i+1}} f\circ \overline{\xi}\}$, the class of $\xi$ vanishes in $\pi_{i}(\widehat{X}_{c})$. Since by Theorem~\ref{th:simpsonpcisoles} the inclusion $\widehat{X}_{d}\subset \widehat{X}_{c}$
induces an isomorphism on $\pi_{i}$, we see that $[\xi]=0$ in $\pi_{i}(\widehat{X}_{d})$. Hence $\pi_{i}(\widehat{X}_{d})=0$.  This proves that the class $b=[Re (\beta )]$ belongs to $\Sigma^{m-1}(\pi_{1}(X))$ (hence to $-\Sigma^{m}(\pi_{1}(X))$ as well).

 The fact that the kernel of $b$ is of type $\mathscr{F}_{m-1}$ if $b$ is rational follows from Theorem~\ref{th:ratfn}. Let us also provide a more direct argument for this result. If $b$ is rational, the image of the integration morphism $\pi_{1}(X)\to \R$ associated to $b$ is cyclic. Consider the associated infinite cyclic covering space $X_0 \to X$ and let $g : X_0 \to \R$ be a primitive of the lift to $X_0$ of the form $Re (\beta)$. Rationality of $b$ implies that the critical values of $g$ are discrete. The map $g$ being proper, each critical level set contains only finitely many critical points. Let $c$ be a regular value of $g$. We can choose an ascending sequence of compact intervals
$$I_0=\left\{c\right\} \subsetneq I_{1}\subsetneq \dots \subsetneq I_j\subsetneq \dots$$
such that $\bigcup_{j\geq 0}I_j = \mathbb{R}$ and $I_j\setminus I_{j-1}$ contains a single critical value of $g$. Since the critical points of $g$ are isolated, Lefschetz theory implies that $g^{-1}(I_j)$ has the homotopy type of $g^{-1}(I_{j-1})$ with finitely many $m$-cells attached to it. Thus, $X_0$ has the homotopy type of a space obtained from $g^{-1}(c)$ by attaching (possibly infinitely many) $m$-cells. Since $g^{-1}(c)$ is a compact manifold and $X_0$ is a $K(\ker(b),1)$, we deduce that $\ker(b)$ is $\mathscr{F}_{m-1}$.

\begin{rem}
 Without assuming $X$ aspherical, the above argument also shows that the inclusion $g^{-1}(c)\hookrightarrow X_0$ is $(m-1)$-connected, i.e. induces an isomorphism on $\pi_i$ for $i<m-1$ and a surjection on $\pi_{m-1}$. The proof of Theorem~\ref{th:simpsonpcisoles} is nothing else than a refinement of this line of argument, taking into account that the cohomology class need not be rational. This means that the set of singular values of $f$ need no longer be discrete. However, the key point is that the set of critical points of $f$ remains discrete in $\widehat{X}$, hence we can still apply Lefschetz theory to reconstruct $f^{-1}([c,\infty))$ from $f^{-1}([d,\infty))$ up to homotopy by attaching cells of dimension $m$. We will explain this in more detail in Section \ref{sec:simpson}.
\end{rem}

To complete the proof of Theorem~\ref{theo:varcomplexabs} we need to show that if $X$ has non-trivial Euler characteristic and $b$ is rational, its kernel is not of type ${\rm FP}_m(\mathbb{Q})$. The following result shows this.
\begin{lemma}\label{lem:notFPm}
Let $Y$ be a $2m$-dimensional closed aspherical real manifold and let $\left[\chi\right]\in S(\pi_1(Y))$ be a character with kernel of type $\mathscr{F}_{m-1}$. Then $b_i^{(2)}(\pi_{1}(Y))=0$ for $2i\neq m$. In particular, Singer's conjecture holds for $Y$. If moreover $b_m^{(2)}(\pi_{1}(Y))\neq 0$, then $\ker(\chi)$ is not of type ${\rm FP}_m(\mathbb{Q})$.
\end{lemma}
\noindent {\it Proof.} Since the kernel of $\chi$ is of type $\mathscr{F}_{m-1}$, the $\ell^2$-Betti numbers $(b_{j}^{(2)}(\pi_{1}(Y)))_{0\le j \le m-1}$ vanish (see~\cite{fisher} or~\cite[Prop. 14]{limp}, this is a consequence of L\"uck's work \cite{luck1998-II}). By Poincar\'e duality, the $\ell^{2}$-Betti numbers $(b_{j}^{(2)}(\pi_{1}(Y)))_{m+1 \le j \le 2m}$ must also vanish. This proves the first two assertions. If $b^{(2)}_{m}(\pi_{1}(Y))\neq 0$ it then follows from another application of \cite[Proposition 14]{limp} that $\ker(\chi)$ is not ${\rm FP}_m(\mathbb{Q})$.\hfill $\Box$

Note that one can provide a direct proof that the kernel of $b$ is not of type ${\rm FP}_m(\mathbb{Q})$ using either complex Morse theory arguments as in \cite{nicolaspy} or homology of cyclic coverings as in \cite[Prop. 21]{limp}. Both arguments provide the stronger result that $H_m(\ker (b),\mathbb{Q})$ is not finite dimensional.

 Before stating our next proposition, we recall that a $G_{\delta}$ is by definition a countable intersection of open sets. 

\begin{prop}\label{prop:lacestod} Let $A$ be a complex torus. Let $U\subset H^{0}(A,\Omega^{1}_{A})$ be the set of holomorphic $1$-forms which do not vanish on any nontrivial subtorus. Then $U$ contains a dense symmetric $G_{\delta}$. Consequently, the set $O=\{ a\in H^{1}(A,\R) | \, a=[Re (\alpha)] \, {\rm with} \, \alpha\in U\}$ contains a dense symmetric $G_{\delta}$ of $H^{1}(A,\R)$. 
\end{prop}

\noindent {\it Proof.} We identify $A$ with $\C^{n}/\Lambda$ where $n={\rm dim}_{\C} \, A$ and $\Lambda < \C^{n}$ is a lattice. The space $H^{0}(A,\Omega^{1}_{A})$ is then identified with the dual space $(\C^{n})^{\ast}$. We define: $$U_{0}=\{\phi \in (\C^{n})^{\ast}\mid \phi (\gamma)\neq 0, \forall \gamma\in \Lambda-\{0\}\}.$$
The set $U_{0}$ is obviously a $G_{\delta}$ and we shall check that $U_0 \subset U$. Let $\phi\in (\C^{n})^{\ast}$ be a holomorphic $1$-form on $A$ which vanishes on a subtorus $T\subset A\cong \C^{n}/\Lambda$ of positive dimension. Then the inverse image of $T$ in $\C^{n}$ is a linear subspace $V$ such that $\phi (V)=0$ and $V\cap \Lambda < V$ is a lattice. This implies that $\phi$ vanishes on a nontrivial element of $\Lambda$, hence $\phi \notin U_{0}$. This concludes the proof.\hfill $\Box$

Proposition~\ref{prop:lacestod} will be used in Section~\ref{sec:defsimplesttype}, to prove Theorem~\ref{bigsigmahd}. Observe that if the torus $A$ in Proposition \ref{prop:lacestod} is a direct product $A=A_1\times \dots \times A_r$ of pairwise nonisogenous simple tori $A_i$, then it contains only finitely many distinct complex subtori. Thus, in this case the set $U$ is open.

We now turn to the proofs of Theorems~\ref{th:fvsigmadeux} and~\ref{th:singerpartial}. 

\noindent {\it Proof of Theorem~\ref{th:fvsigmadeux}.} Let $X$ be as in the statement of the theorem. For $m=1$, $X$ is a closed Riemann surface and the assertion is well-known. So we assume that $m\ge 2$ and, by contradiction, that $\Sigma^{m}(\pi_{1}(X))$ is nonempty. Since this set is open, Propositions~\ref{prop:alsago} and~\ref{prop:lacestod} imply that there exists a holomorphic $1$-form with finitely many zeros on $X$ such that the cohomology class $[Re (\alpha) ]$ lies in $\Sigma^{m}(\pi_{1}(X))$. We let $\chi : \pi_{1}(X)\to \R$ be the induced character. 

To make our argument more transparent, we shall consider a perturbation $\alpha'$ of $\alpha$ constructed as follows. For each point  $p\in X$ such that $\alpha_p =0$ we pick two small balls $V_{p} \subset \overline{V_{p}}\subset U_{p}$ centered at $p$ (in some chart) such that the $U_p$'s are pairwise disjoint. By adding a small generic complex linear form to $\alpha$ in $U_p$ we obtain a perturbation with only nondegenerate singularities. This form on $U_{p}$ can be perturbed in $U_p -V_p$ in a $C^{\infty}$ way to coincide with $\alpha$ near the boundary of $U_p$. This allows us to construct a form $\alpha'$ close to $\alpha$ (and cohomologous to it), which is holomorphic everywhere except in the set 
$$\bigcup_{p, \alpha_{p}=0}U_{p}-V_{p},$$
and whose zeros are nondegenerate and contained in the union of the $V_p$'s. Theorem~\ref{th:simpsonpcisoles} still applies to $\alpha'$. Indeed the proof only uses that the form is holomorphic near its critical points (see Section~\ref{sec:simpson}). Since $\alpha$ and $\alpha '$ define the same cohomology class, we have $[Re ( \alpha ') ]\in \Sigma^{m}(\pi_{1}(X))$. 

Let now $\widehat{X}$ be the universal cover of $X$ and let $f : \widehat{X}\to \R$ be a primitive of the lift of the form $Re (\alpha ')$ to $\widehat{X}$. As before, for a real number $d$, we write
$$\widehat{X}_{d}:=f^{-1}([d,\infty)).$$
We first claim that the $(m-1)$-th homotopy group of $\widehat{X}_{d}$ vanishes for every real number $d$. By the definition of the $m$-th BNSR invariant, there exists a real number $r\ge 0$ such that the map
\begin{equation}\label{maphomotopygroups}  
\pi_{i}(\widehat{X}_{0})\to \pi_{i}(\widehat{X}_{-r})
\end{equation}
has a one point image for $i\le m-1$. Let $\xi \in \pi_{m-1}(\widehat{X}_{-r})$. According to Theorem~\ref{th:simpsonpcisoles}, $\xi$ lies in the image of the map~\eqref{maphomotopygroups}. Since this map is trivial, $\xi=0$ and thus $\pi_{m-1}(\widehat{X}_{-r})$ is trivial. It follows that for every element $g\in \pi_{1}(X)$, $g(\widehat{X}_{-r})=\widehat{X}_{\chi (g)-r}$ also has trivial $\pi_{m-1}$. Given an arbitrary real number $d$, there is an element $g\in \pi_{1}(X)$ such that $\chi (g)-r> d$. Theorem~\ref{th:simpsonpcisoles} implies that the inclusion 
$$\widehat{X}_{\chi (g)-r} \subset \widehat{X}_{d}$$
induces a surjection on $\pi_{m-1}$, hence $\pi_{m-1}(\widehat{X}_{d})=0$. This proves the claim.

 We will now finish the proof by an argument which is very similar in spirit to arguments appearing in~\cite{kapovich,nicolaspy}. Since the Euler characteristic of $X$ is nonzero, we can pick two real numbers $c>d$ such that there exists at least one zero of $\alpha '$ in the open set $\{ d < f < c\}$. The proof of Theorem~\ref{th:simpsonpcisoles} shows that $\widehat{X}_{d}$ has the homotopy type of a space $W_{c}$ obtained from $\widehat{X}_{c}$ by gluing some $m$-dimensional balls $(B_{i})_{i\in I}$ along a non-empty family of disjoint spheres $(S_{i})_{i\in I}$ contained in the set $\{ f=c\}$. We pick one of these spheres, say $S_{i_{0}}$, and consider the corresponding ball $B_{i_{0}}$. Since $\pi_{m-1}(\widehat{X}_{c})$ is trivial, there exists a map $v: B^{m} \to \widehat{X}_{c}$ such that the restriction of $v$ to the boundary of $B^m$ is a homeomorphism onto $S_{i_{0}}$. By ``gluing" $v$ and a parametrization of the ball $B_{i_{0}}$, we obtain a map $$v^{1} : S^{m}\to W_{c}$$ from the $m$-dimensional sphere to the space $W_c$. Writing $W_{c}$ as the union of $\widehat{X}_{c}$ and some $m$-dimensional balls and using the Mayer-Vietoris exact sequence we see that $v^{1}$ defines a nontrivial homology class $v^{1}_{\ast}([S^{m}])$ in $H_{m}(W_c,\Q)$. Applying the proof of Theorem~\ref{th:simpsonpcisoles} again, we see that $\widehat{X}$ is obtained from $\widehat{X}_{d}$ by gluing successively some $m$-dimensional balls. In particular the map
 $$H_{m}(\widehat{X}_{d},\Q)\to H_{m}(\widehat{X},\Q)$$
 is injective. Since $\widehat{X}_d$ has the same homotopy type as $W_c$, the existence of the nontrivial element $v^{1}_{\ast}([S^{m}])$ of $H_{m}(W_c,\Q)$ thus contradicts the asphericity of $\widehat{X}$. This concludes the proof.\hfill $\Box$

As a side remark, we state the following proposition, which follows easily from Gromov's theorem characterizing K\"ahler groups with nonzero first $\ell^{2}$-Betti number~\cite{abckt,gromovcras}. 

\begin{prop}\label{prop:asphbeun} Let $X$ be a closed aspherical K\"ahler manifold with ${\rm dim}_{\C}\, X\ge 2$. Then the first $\ell^{2}$-Betti number of $\pi_{1}(X)$ is zero. 
\end{prop}

\noindent {\it Proof.} We assume by contradiction that $X$ has positive first $\ell^{2}$-Betti number. A result of Gromov~\cite{gromovcras} then implies that there exists a finite covering space $X_1 \to X$ and a holomorphic map with connected fibers $$p : X_{1}\to \Sigma$$ onto a closed hyperbolic Riemann surface such that the fundamental group of every smooth fiber of $p$ has finite image in $\pi_{1}(X_{1})$. We refer the reader to~\cite[Ch. 4]{abckt} for a proof. Let $F$ be a smooth fiber of $p$. Since the image
$${\rm Im}(\pi_{1}(F)\to \pi_{1}(X_{1}))$$
is finite, there exists a finite covering space $F_{1}\to F$ such that $F_1$ can be lifted to the universal cover $\widehat{X}$ of $X$. Since $X$ is K\"ahler, the image of this lift defines a nontrivial homology class in $\widehat{X}$. This contradicts the asphericity of $X$ and finishes the proof.\hfill$\Box$

As a consequence of this proposition we have:

\begin{cor}
Let $X$ be a closed aspherical K\"ahler surface. Then $X$ satisfies Singer's conjecture.    
\end{cor} 

\noindent {\it Proof.} According to Proposition~\ref{prop:asphbeun}, we have $b_{1}^{(2)}(\pi_{1}(X))=0$. By Poincar\'e duality, $b_{3}^{(2)}(\pi_{1}(X))=0$, hence $X$ satisfies Singer's conjecture.\hfill $\Box$

Using Theorem~\ref{th:ratfn} and Lemma~\ref{lem:notFPm} one sees easily that if $X$ is as in the corollary and has nonzero Euler characteristic, the intersection $\Sigma^{2}(\pi_{1}(X))\cap -\Sigma^{2}(\pi_{1}(X))$ is empty. This was already observed by different means by Friedl and Vidussi~\cite{friedlvidussi}, see Proposition 3.4 there and the remark following it. However one has actually that 
$$\Sigma^{k}(\pi_{1}(M))\cap -\Sigma^{k}(\pi_{1}(M))=\emptyset$$
for any aspherical closed $2k$-manifold $M$ of nonzero Euler characteristic. This follows easily from Milnor's work~\cite{milnor} and from~\cite[Prop. 21]{limp}. 

\noindent {\it Proof of Theorem~\ref{th:singerpartial}.} Let $X$ and $\alpha$ be as in the statement of the theorem. According to Theorem~\ref{theo:varcomplexabs}, the cohomology class $a$ of the real $1$-form $Re(\alpha)$ lies in $\Sigma^{m-1}(\pi_{1}(X))\cap -\Sigma^{m-1}(\pi_{1}(X))$. In particular $\Sigma^{m-1}(\pi_{1}(X))\cap - \Sigma^{m-1}(\pi_{1}(X))$ is non-empty. Since this set is open, we can pick a rational class $a_{0}$ in it. According to Theorem~\ref{th:ratfn} the kernel of $a_{0}$ is then of type $\mathscr{F}_{m-1}$ and Lemma \ref{lem:notFPm} implies Singer's conjecture for $X$. \hfill $\Box$

\subsection{Simpson's theorem}\label{sec:simpson}
For the reader's convenience, we include a brief account of Simpson's work~\cite{simpson} in this section, providing a proof of Theorem \ref{th:simpsonpcisoles}. We restrict ourselves to the situation required in our work, although Simpson's results are stated in greater generality. Moreover, we assume that $\beta$ has only nondegenerate singularities. The general case follows by a simple perturbation argument, see e.g.~\cite[\S 2]{nicolaspy}. 

Let $\{p_{1}, \ldots , p_{r}\}$ be the set of zeros of $\beta$. We choose disjoint open sets $(O_{i})_{1\le i\le r}$ in $Y$ such that $p_{i}\in O_{i}$ and such that for each $i\in \{1, \ldots r\}$, there exists a biholomorphic map 
$$\phi_{i} : O_{i}\to B(0,2)\subset \C^m$$
which takes $p_i$ to $0$. We denote by $h_{i} : O_{i} \to \C$ the primitive of $\beta $ on $O_{i}$ such that $h_{i}(p_{i})=0$ and let $f_{i}=h_{i}\circ \phi_{i}^{-1}$. Standard arguments from the study of Milnor fibrations show that there exists $\varepsilon >0$ such that for each $i\in \{1, \ldots , r\}$ the restriction of $f_i$ to the boundary of the ball $B(0,1)$ is a submersion at each point of the set $\{\vert Re f_{i}\vert \le \varepsilon\}$. We now state two lemmas whose proofs are left to the reader. The first one follows from standard arguments from Morse or Lefschetz theory~\cite[Ch. 14]{voisin}; the second follows easily from the fact that $f_i |_{\partial B(0,1)}$ is a submersion along $\partial B(0,1)\cap \{\vert Re f_{i}\vert \le \varepsilon\}$. From now on we fix a positive number $\delta \le \varepsilon$. 

If $N$ is a closed manifold with boundary endowed with a submersion $p : N\to I$ where $I=[s_{1},s_{2}]$ is a closed interval of $\R$, each trivialization $\Phi : N\to p^{-1}(s_{2}) \times I$ of $p$ provides a canonical retraction by deformation of $N$ onto $p^{-1}(s_{2})$ given (in the trivialization) by the map $(x,s)\mapsto x$. We shall say that this retraction is induced by the trivialization $\Phi$. 

\begin{lemma}\label{lemma:ulcl} For all real numbers $\mu$ and $\lambda$ such that $-\delta \le \mu \le 0 < \lambda \le \delta$ the set $B(0,1)\cap \{\mu \le Re (f_{i})\le \lambda \}$ deformation retracts onto the union of the level set $B(0,1)\cap \{Re(f_{i})=\lambda\}$ and a ball of dimension $m$ glued along an $(m-1)$-dimensional sphere contained in $B(0,1)\cap \{Re(f_{i})=\lambda\}$. 

This retraction can be chosen to coincide on $\partial B(0,1)\cap \{\mu \le Re(f_{i})\le \lambda\}$ with the retraction induced by any trivialization of the bundle $\partial B(0,1)\cap \{\mu \le Re(f_{i})\le \lambda\}\to [\mu , \lambda]$.  
\end{lemma}

In the lemma below we write for $-\delta \le u \le u' \le \delta$, 
\begin{equation}\label{eq:defouverts}
U_{i,u,u'}=\phi_{i}^{-1}(\overline{B(0,1)})\cap \{u\le  Re(h_{i}) \le u'\},
\end{equation}
and denote by ${\rm Int}(U_{i,u,u'})$ the interior of $U_{i,u,u'}$. 

\begin{lemma}\label{lemma:vfield} There exists a smooth vector field $V$ on $Y$ such that $0\le Re (\beta )(V)\le 1$ on $Y$, $Re (\beta) (V)=1$ outside of 
$$\cup_{i=1}^{r}{\rm Int}(U_{i,-\frac{\delta}{10},\frac{\delta}{10}})$$ and such that $V$ is tangent to $\phi_{i}^{-1}(\partial B(0,1))$ along $\phi_{i}^{-1}(\partial B(0,1)) \cap \{\vert Re (f_{i})\vert \le \delta\}$. 
\end{lemma}
 
Let $\pi : \widehat{Y} \to Y$ be the projection, $h : \widehat{Y} \to \C$ be a primitive of the lift of $\beta$ to $\widehat{Y}$ and let $f=Re(h)$. Let $\widehat{V}$ be the lift to $\widehat{Y}$ of the vector field $V$ from Lemma~\ref{lemma:vfield}. If $x\in \widehat{Y}$ is a critical point of $h$, the image of $x$ in $Y$ equals one of the zeros of $\beta$, say $p_{i}$. If $-\delta \le u \le u'\le \delta$, we let $U_{x,u,u'}$ be the component of the preimage by $\pi$ of $U_{i,u,u'}$ containing $x$. Since $\pi$ identifies $U_{x,u,u'}$ and $U_{i,u,u'}$ as well as $h$ and $h_{i}$ (up to translation) we can apply Lemma~\ref{lemma:ulcl} to the map $h : U_{x,u,u'}\to \C$. Note that using our convention~\eqref{eq:defouverts}, we have $u\le f-f(x)\le u'$ on $U_{x,u,u'}$ since $f$ and $Re(h_{i})\circ \pi$ differ by the constant $f(x)$ on $U_{x,u,u'}$. Let $c$ be a real number and let $\Lambda_{c}$ be the set of critical points $x$ of $h$ such that $$[c,c+\frac{\delta}{10}]\cap [f(x)-\frac{\delta}{10},f(x)+\frac{\delta}{10}]$$ is nonempty. We now define a continuous map $$F : \{f\ge c \} \to \{ f\ge c\}$$ whose image will be the union of $\{f\ge c+\frac{\delta}{10}\}$ together with countably many $m$-dimensional cells glued to $\{f\ge c+\frac{\delta}{10}\}$ along their boundary. The map $F$ will be a retraction by deformation onto its image. This implies the conclusion of Theorem~\ref{th:simpsonpcisoles} whenever $c\le d\le \frac{\varepsilon}{10}$. The general case follows by applying this step finitely many times. Note that $\varepsilon$ is fixed once and for all and only depends on $Y$ and $\beta$.   

The map $F$ is built as follows. If a point $z$ does not belong to the set 
\begin{equation}\label{eq:ouvertopp}
\bigcup_{x\in \Lambda_{c}}U_{x,c-f(x),c-f(x)+\frac{\delta}{10}},
\end{equation}
one follows the flow line of $\widehat{V}$ starting from $z$ until one reaches a point of the level set $f=c+\frac{\delta}{10}$. This defines $F(z)$. If $z$ belongs to the set~\eqref{eq:ouvertopp}, one constructs $F$ as follows, noting that by Lemma \ref{lemma:vfield} we have the necessary freedom to define it so that it is continuous on $\left\{f\ge c\right\}$. If $z\in U_{x,c-f(x),c-f(x)+\frac{\delta}{10}}$ and $c\le f(x) < c+\frac{\delta}{10}$, one applies Lemma~\ref{lemma:ulcl} to build a retraction of $U_{x,c-f(x),c-f(x)+\frac{\delta}{10}}$ onto the union of $U_{x,c-f(x),c-f(x)+\frac{\delta}{10}} \cap \{ f=c+\frac{\delta}{10}\}$ with an $m$-dimensional sphere. If $f(x)=c+\frac{\delta}{10}$, one can retract $U_{x,c-f(x),c-f(x)+\frac{\delta}{10}}$ onto $U_{x,c-f(x),c-f(x)+\frac{\delta}{10}}\cap \{ f=c+\frac{\delta}{10}\}$. The case where $f(x)\notin [c, c+\frac{\delta}{10}]$ is even simpler since the map $U_{x,c-f(x),c-f(x)+\frac{\delta}{10}} \to [c,c+\frac{\delta}{10}]$ is a locally trivial fibration in this case. This completes the proof of Theorem~\ref{th:simpsonpcisoles}.


\section{Arithmetic lattices of the simplest type and Albanese maps}\label{sec:chlattices}

We start this section by defining the Albanese map of a compact K\"ahler manifold and studying it for arithmetic quotients of complex balls in Section~\ref{sec:virtuallyimm}. After that we turn to the definition of arithmetic lattices of the simplest type and to the proof of Theorem~\ref{bigsigmahd} and Corollaries~\ref{corollairerationalclass} and~\ref{cor:sigmadeuxvide} in Section~\ref{sec:defsimplesttype}. 

\subsection{Albanese maps are virtually immersions}\label{sec:virtuallyimm}

Let $X$ be a closed K\"ahler manifold and let $H^{0}(X,\Omega_{X}^{1})$ be the space of holomorphic $1$-forms on $X$. If $\theta : [0,1] \to X$ is a path in $X$ we denote by $i(\theta) \in (H^{0}(X,\Omega_{X}^{1}))^{\ast}$ the linear map on $H^{0}(X,\Omega_{X}^{1})$ taking a form $\alpha$ to the integral
$$\int_{\theta}\alpha.$$
Since holomorphic forms on $X$ are closed, this only depends on the homotopy class of $\theta$ relative to its endpoints. Thus, one can also define $i(u)$ for $u\in H_{1}(X,\Z)$. The kernel of the map 
$$i : H_{1}(X,\Z)\to (H^{0}(X,\Omega_{X}^{1}))^{\ast}$$
is the torsion subgroup of $H_{1}(X,\Z)$ and its image is a lattice in $(H^{0}(X,\Omega_{X}^{1}))^{\ast}$. The Albanese torus of $X$ is defined as
\begin{equation}
 {\rm Alb}(X):=(H^{0}(X,\Omega_{X}^{1}))^{\ast}/i(H_{1}(X,\Z)).
\end{equation}
For a fixed point $x_0 \in X$ we define a map 
$${\rm alb}_{X} : X \to {\rm Alb}(X)$$
by setting ${\rm alb}_{X}(x)= i(\theta_{x}) \, {\rm mod}~ i(H_{1}(X,\Z))$, where $\theta_{x}$ is any continuous path going from $x_0$ to $x$. This does not depend on the choice of $\theta_{x}$. The resulting map ${\rm alb}_{X}$ is holomorphic. We refer the reader to~\cite[\S 12.1.3]{voisin} for more details on these notions. By construction the differential of the map ${\rm alb}_{X}$ at a point $x\in X$ is the evaluation map 
$$\begin{array}{ccc}
T_{x}X & \to & (H^{0}(X,\Omega_{X}^{1}))^{\ast} \\
v & \mapsto & (\alpha \mapsto \alpha_{x}(v)).\\ 
\end{array}$$
Hence we have:
\begin{lemma}\label{lemma:imm}
Let $x\in X$. The linear map $d {\rm alb}_{X} (x)$ is injective if and only if the evaluation map 
$$\begin{array}{ccc}
H^{0}(X,\Omega_{X}^{1}) & \to & (T_{x}X)^{\ast} \\
\alpha & \mapsto & \alpha_{x} \\ 
\end{array}$$
is onto. 
\end{lemma}

We now turn to the study of the Albanese map in the case of a quotient of the unit ball $B \subset \mathbb{C}^{m}$ by an arithmetic lattice. Let $\Gamma < {\rm PU}(m,1)$ be a cocompact torsionfree lattice. We recall that the {\it commensurator} of $\Gamma$ is defined as follows:
$${\rm Comm}(\Gamma)=\{g\in {\rm PU}(m,1) \mid \Gamma \cap g\Gamma g^{-1} {\rm has} \,\, {\rm finite} \,\, {\rm index} \,\, {\rm in} \,\, {\rm both} \,\, \Gamma \,\,{\rm and} \,\, g\Gamma g^{-1}\}.$$
This is a subgroup of $\Gamma$ and a well-known theorem of Margulis states that $\Gamma$ is arithmetic if and only if ${\rm Comm}(\Gamma)$ is a dense subgroup of ${\rm PU}(m,1)$. Below we will make the assumption that $\Gamma$ is arithmetic but this hypothesis will only play a role through the density of the commensurator of $\Gamma$ so that the reader can take the conclusion of Margulis' theorem as a definition of arithmeticity if they so desire.

We assume that $b_{1}(\Gamma)>0$ and that $\Gamma$ is arithmetic. In that case it is well-known that the virtual first Betti number of $\Gamma$ must be infinite~\cite{agol0,venk}. The following theorem relies on similar ideas and gives a geometric application. 

\begin{main}\label{th:immersion} (Eyssidieux) Assume that $\Gamma$ is arithmetic and has positive first Betti number. Then there exists a finite index subgroup $\Gamma_{0}<\Gamma$ with the property that the Albanese map of $X_{\Gamma_{0}}$ is an immersion. 
\end{main}

This result appears in~\cite{eyssidieux}. More precisely, in that work the author proves Theorem~\ref{th:immersion} for non-uniform lattices and gives applications of this result to the Shafarevich conjecture of holomorphic convexity for certain toroidal compactifications. The cocompact case is analogous. Since the proof is quite simple, we include it here for the reader's convenience. We also note that part of Eyssidieux's arguments are identical to the ones in~\cite{agol0}.  

In  what follows we shall say that two lattices $\Gamma_{1}, \Gamma_{2} < {\rm PU}(m,1)$ are commensurable if $\Gamma_{1}\cap \Gamma_{2}$ has finite index in both $\Gamma_{1}$ and $\Gamma_{2}$. When using the abstract notion of commensurability from Definition~\ref{def:cor}, we shall say it explicitly. 

We fix $\Gamma < {\rm PU}(m,1)$ as in Theorem~\ref{th:immersion}. Let $\Omega_{B}^{1}$ be the space of all holomorphic $1$-forms on $B$. For a lattice $\Lambda < {\rm PU}(m,1)$, let $\Omega_{B,\Lambda}^{1}\subset \Omega_{B}^{1}$ be the subspace of $\Lambda$-invariant forms. We define a subset $\mathscr{L}\subset \Omega_{B}^{1}$ as follows. A holomorphic $1$-form $\alpha$ on $B$ belongs to $\mathscr{L}$ if and only if there exists a cocompact lattice $\Lambda < {\rm PU}(m,1)$ such that $\Lambda$ is commensurable to $\Gamma$ and $\alpha$ is invariant under $\Lambda$. Following Agol~\cite{agol0} and Eyssidieux~\cite{eyssidieux}, we consider the linear subspace
$$V_{0}\subset \Omega^1_B$$
spanned by $\mathscr{L}$ and let $V$ be the closure of $V_0$ for the topology of uniform convergence on compact sets. We first observe that the set $\mathscr{L}$ is invariant under the action of ${\rm Comm}(\Gamma)$. Indeed, if $\alpha \in \mathscr{L}$ is invariant under a lattice $\Lambda$ commensurable with $\Gamma$, and if $g\in {\rm Comm}(\Gamma)$, then $(g^{-1})^{\ast}\alpha$ is invariant under the lattice $g\Lambda g^{-1}$ which is still commensurable with $\Gamma$. Hence $g(\mathscr{L})=\mathscr{L}$ and consequently $g(V_{0})=V_{0}$ and $g(V)=V$ for all $g\in {\rm Comm}(\Gamma)$. The density of ${\rm Comm}(\Gamma)$ in ${\rm PU}(m,1)$ then implies that the space $V$ is ${\rm PU}(m,1)$-invariant. 

Since the intersection of finitely many lattices commensurable with $\Gamma$ is again a lattice commensurable with $\Gamma$, the following lemma is clear. 

\begin{lemma}\label{lemma:invsubgroup} Let $W\subset V_0$ be a finite dimensional vector subspace. Then there exists a lattice $\Lambda < {\rm PU}(m,1)$ commensurable with $\Gamma$ such that $\gamma^{\ast}\alpha=\alpha$ for all $\alpha \in W$ and $\gamma \in \Lambda$. 
\end{lemma}

\begin{lemma}\label{lemma:surjclv} Let $p\in B$. Let $ev_{p} : V \to (\C^{m})^{\ast}$ be the evaluation map at $p$. Then $ev_p$ is onto.
\end{lemma}

\noindent {\it Proof of Lemma~\ref{lemma:surjclv}.} Since $V$ is ${\rm PU}(m,1)$-invariant, it is enough to prove the Lemma for $p=o$, the origin of the ball. If the image $ev_{o}(V)\subset (\mathbb{C}^{m})^{\ast}$ is equal to $0$, then $ev_{p}$ would be equal to $0$ for each point $p\in B$ and $V$ would be reduced to zero. This is impossible since $b_{1}(\Gamma) >0$. Hence the image of $ev_{o}$ is nonzero. Since $ev_{o}(V)\subset (\mathbb{C}^{m})^{\ast}$ is invariant under the natural action of ${\rm U}(m)$, we must have $ev_{o}(V)=(\mathbb{C}^{m})^{\ast}$. This concludes the proof.\hfill $\Box$

For a lattice $\Lambda < {\rm PU}(m,1)$, we now define the following subset of the ball:
\begin{equation}
Z_{\Lambda}=\{x\in B \mid d{\rm alb}_{X_{\Lambda}} \,\, {\rm is} \,\, {\rm not} \,\, {\rm injective} \,\, {\rm at} \,\, x \,\, {\rm mod}~\Lambda\}.
\end{equation}
We make the following observations. Given a point $p\in B$, Lemma~\ref{lemma:surjclv} implies that there exist elements $\alpha_{1}, \ldots , \alpha_{m}$ in $V$ such that the evaluations $(\alpha_{i}(p))_{1\le i \le m}$ generate $(\mathbb{C}^{m})^{\ast}$. Since this is an open condition, we can actually assume that $\alpha_{1}, \ldots , \alpha_{m}$ belong to $V_{0}$. Note that the evaluations $\alpha_{1}(q), \ldots , \alpha_{m}(q)$ will then be linearly independent for $q$ in a neighborhood of $p$. This implies that for each compact subset $M\subset B$, there exists a finite dimensional subspace $W\subset V_{0}$ such that $ev_{p}(W)=(\mathbb{C}^{m})^{\ast}$ for each point $p\in M$. According to Lemma~\ref{lemma:invsubgroup}, there exists a lattice $\Lambda$ commensurable with $\Gamma$ such that $W\subset \Omega_{B,\Lambda}^{1}$. Thanks to Lemma~\ref{lemma:imm}, this implies that 
\begin{equation}
M\cap Z_{\Lambda}=\emptyset.
\end{equation} 
In other words, we have proved: 

\begin{prop}\label{prop:keyimmersion} For each compact subset $M\subset B$, there exists a lattice $\Lambda < {\rm PU}(m,1)$ commensurable with $\Gamma$ such that $Z_{\Lambda}$ does not intersect $M$. 
\end{prop}

We now pick a compact fundamental domain $K\subset B$ for the action of $\Gamma$ on $B$. According to Proposition~\ref{prop:keyimmersion}, there exists  a lattice $\Lambda < {\rm PU}(m,1)$ commensurable with $\Gamma$ and such that $Z_{\Lambda} \cap K=\emptyset$. Let $\Gamma_{0} \lhd \Gamma$ be a normal finite index subgroup such that $\Gamma_{0} \subset \Lambda \cap \Gamma$. Since $\Gamma_0 < \Lambda$ we must have
$$Z_{\Gamma_{0}} \subset Z_{\Lambda}.$$
Hence $Z_{\Gamma_{0}}$ does not intersect $K$. But $Z_{\Gamma_{0}}$ is $\Gamma$-invariant since $\Gamma_{0}\lhd \Gamma$. If $x\in Z_{\Gamma_{0}}$, there exists $\gamma \in \Gamma$ such that $\gamma \cdot x\in K$. But $\gamma \cdot x$ also lies in $Z_{\Gamma_{0}}$. We thus obtain a contradiction. This shows that $Z_{\Gamma_{0}}=\emptyset$, implying that the Albanese map of $X_{\Gamma_{0}}$ is an immersion. This concludes the proof of Theorem~\ref{th:immersion}. 

\begin{rem} We shall describe below a class of arithmetic lattices in ${\rm PU}(m,1)$ which are known to have positive virtual first Betti number. Note though that there are other classes of arithmetic lattices for which we currently don't know if they have positive virtual first Betti number. This includes a class for which the first Betti number is known to vanish on all congruence subgroups, while the existence of noncongruence finite index subgroups is open~\cite{clozel}. We emphasize that all the results from this paper which are stated for arithmetic lattices of the simplest type in fact hold for all arithmetic lattices in ${\rm PU}(m,1)$ with positive first Betti number. Indeed our proofs only use the arithmeticity, not the specific arithmetic construction.
\end{rem}

\subsection{Conclusion of the proofs}\label{sec:defsimplesttype} 

In this section we prove Theorem~\ref{bigsigmahd} and Corollaries~\ref{corollairerationalclass} and~\ref{cor:sigmadeuxvide}. We start by recalling the definition of arithmetic lattices of the simplest type in the group ${\rm PU}(m,1)$. Let $F\subset \R$ be a totally real number field and $E\subset \C$ be a purely imaginary quadratic extension of $F$. Let $V=E^{m+1}$ and let $H : V \times V \to E$ be a hermitian form. We assume that the extension of $H$ to $V\otimes \C$ has signature $(m,1)$ and that for every embedding $\sigma : E \to \C$ with $\sigma|_{F}$ distinct from the identity of $F$, the twisted Hermitian form $H^{\sigma}$ has signature $(m+1,0)$ over $\C$. Let $\mathscr{O}_{E}$ be the ring of integers of $E$ and $U(H, \mathscr{O}_{E})$ be the group of $(m+1)\times (m+1)$ matrices with coefficients in $\mathscr{O}_{E}$ which preserve the Hermitian form $H$. The group $U(H,\mathscr{O}_{E})$ is a lattice in the group $U(V\otimes \C , H)$ of automorphisms of the space $(V\otimes \C,H)$, which is cocompact if and only if $F$ is distinct from $\Q$. It is these lattices and the ones commensurable to them that are usually called ``of the simplest type". See~\cite[\S VIII.5]{bowa} for more details. 

A result of Kazhdan~\cite{kazhdan} states that when $F\neq \mathbb{Q}$, the group $U(H,\mathscr{O}_{E})$ always admits congruence subgroups with positive first Betti number. This result is also exposed in~\cite{bowa} and has been extended by Shimura to the noncocompact case~\cite{shimura}.    

\noindent {\it Proof of Theorem~\ref{bigsigmahd}.} Let $\Gamma < {\rm PU}(m,1)$ be a torsion-free cocompact arithmetic lattice of the simplest type. Applying Kazhdan's result and then Theorem~\ref{th:immersion} we obtain a finite index subgroup $\Gamma_{0}< \Gamma$ such that the Albanese map of the manifold
$$X_{\Gamma_{0}}=B/\Gamma_{0}$$
is an immersion. Fix a finite index subgroup $\Gamma_{1}< \Gamma_{0}$. The Albanese map of $X_{\Gamma_{1}}$ is also an immersion since there is a natural commutative square
$$\xymatrix{X_{\Gamma_{1}} \ar[d] \ar[r] & {\rm Alb}(X_{\Gamma_{1}}) \ar[d] \\
X_{\Gamma_{0}} \ar[r] & {\rm Alb}(X_{\Gamma_{0}}). \\}$$
Theorem~\ref{theo:varcomplexabs} implies that if $\alpha$ is a holomorphic $1$-form on ${\rm Alb}(X_{\Gamma_{1}})$ which does not vanish on any nontrivial subtorus, the class
$$a=[Re ({\rm alb}_{X_{\Gamma_{1}}}^{\ast} \alpha)]$$
lies in $\Sigma^{m-1}(\Gamma_{1})$. Since the set of such classes is dense by Proposition~\ref{prop:lacestod}, this concludes the proof of the first statement of the theorem. The affirmation about rational classes follows from Theorem~\ref{th:ratfn} and Lemma~\ref{lem:notFPm}, together with the fact that $B/\Gamma_{1}$ has nonzero Euler characteristic. Again, we could also have appealed to~\cite[Prop. 21]{limp}.\hfill $\Box$

\noindent {\it Proof of Corollary~\ref{corollairerationalclass}.} The following arguments are classical. We will show that arithmetic lattices of the simplest type in ${\rm PU}(m,1)$ form infinitely many commensurability classes in the sense of Definition~\ref{def:cor}. We recall that the {\it adjoint trace field} of a lattice $\Gamma < {\rm PU}(m,1)$ is the field generated by the traces of the transformations 
$${\rm Ad}(\gamma) : {\rm Lie}({\rm PU}(m,1)) \to {\rm Lie}({\rm PU}(m,1))\;\;\;\; (\gamma \in \Gamma).$$
We now observe that if $\Gamma_1$ and $\Gamma_2$ are two lattices in the group ${\rm PU}(m,1)$ ($m\ge 2$) which are commensurable in the sense of Definition~\ref{def:cor}, Mostow's rigidity theorem implies that there exists an element $g\in {\rm PU}(m,1)$ such that the intersection
$$g\Gamma_{1}g^{-1}\cap \Gamma_{2}$$
has finite index in both $g\Gamma_{1}g^{-1}$ and $\Gamma_{2}$. This implies that $\Gamma_{1}$ and $\Gamma_{2}$ have the same adjoint trace field~\cite[Prop. 12.2.1]{dm}. But the adjoint trace field of the lattice $U(H,\mathscr{O}_{E})$ built above is known to be the field $F$~\cite[12.2.5]{dm}. Since there are infinitely many possibilites for $F$, this proves the corollary.\hfill $\Box$ 

\noindent {\it Proof of Corollary~\ref{cor:sigmadeuxvide}.} Let $\Gamma_1 < \Gamma$ be a finite index subgroup such that the Albanese map of the manifold $X_{\Gamma_{1}}$ is an immersion. We denote by $i : \Gamma_1 \to \Gamma$ the inclusion. If $\chi\in H^{1}(\Gamma,\R)-\{0\}$ is such that $[\chi]\in \Sigma^{m}(\Gamma)$, then $[\chi \circ i]\in \Sigma^{m}(\Gamma_{1})$; this follows directly from the definition of $\Sigma^{m}$. But $\Sigma^{m}(\Gamma_{1})$ is empty by Theorem~\ref{th:fvsigmadeux}. This forces $\Sigma^{m}(\Gamma)$ to be empty and finishes the proof.\hfill $\Box$

\begin{rem}
Let $\Gamma < {\rm PU}(m,1)$ be a cocompact arithmetic lattice of the simplest type. Theorem~\ref{bigsigmahd} improves on Stover's work~\cite{stover} who proved that the invariant $\Sigma^{1}$ of deep enough finite index subgroups of $\Gamma$ is nonempty. In particular, for $m=2$, Theorem~\ref{bigsigmahd} is a consequence of~\cite{delzant,stover}. Indeed, Delzant's work shows that the invariant $\Sigma^{1}$ of a K\"ahler group, when nonempty, is the complement of the union of finitely many proper subspheres of the character sphere.
\end{rem}

\bigskip
\bigskip
\begin{tiny}
\begin{tabular}{lllllll}
Claudio Llosa Isenrich & & & & Pierre Py \\
Faculty of Mathematics & & & & IRMA \\
Karlsruhe Institute of Technology & & & & Universit\'e de Strasbourg \& CNRS \\
76131 Karlsruhe, Germany  & & & & 67084 Strasbourg, France \\
claudio.llosa at kit.edu & & & & ppy at math.unistra.fr \\    
\end{tabular}

\end{tiny}

\end{document}